\documentclass[reqno,fleqn]{amsart}
\usepackage{amsfonts,amsthm,amssymb}
\usepackage{lipsum}
\usepackage[utf8]{inputenc}
\usepackage[english]{babel}

\usepackage{lmodern,textcomp}
\usepackage{amssymb}
\usepackage{amsthm}
\usepackage{amsmath}
\usepackage{mathtools}
\usepackage{latexsym}
\usepackage{tikz}
\usepackage{fancyvrb}
\usepackage{epsfig}
\usepackage{graphicx}
\usepackage{amsmath}
\usepackage{latexsym}
\usepackage{tikz}
\usepackage{rotating}
\usetikzlibrary{positioning}
\usepackage{subcaption}
\usepackage{datetime}
\newtheorem{theorem}{Theorem}[section]

\newtheorem{cor}[theorem]{Corollary}
\newtheorem{definition}[theorem]{Definition}
\newtheorem{conj}{Conjecture}
\newtheorem{exm}{Example}

\newtheorem{note}[theorem]{Note}

\title[Cylindrical Grid Graphs $P_m \Box C_n$ are Non-Distance Magic]{Cylindrical Grid Graphs $P_m \Box C_n$ are Non-Distance Magic}

\author{\sc Sajidha P, V. Vilfred Kamalappan and Julia K. Abraham} 
\address{Department of Mathematics, ~Central ~University ~of~ Kerala,~Periye \\  Kasaragod,~ Kerala, ~India~ - ~671 316.}

\email{sajinada555@gamil.com}
\email{vilfredkamalv@cukerala.ac.in}
\email{juliakabraham92@gmail.com}


\date{}

\begin{document}
	\begin{abstract} A bijective mapping $f: V(G) \rightarrow \left\{1,2,\ldots,n\right\}$ is called a \emph{Distance Magic Labeling (DML) of $G$} if ~ ${\sum_{v \in N(u)}} f(v) $ is a constant for all  $u\in V(G)$ where $G$ is a simple graph of order $n$ and $N(u)$ = $\{v\in V(G):$ $uv\in E(G)\}$. Graph $G$ is called a \emph{Distance Magic Graph (DMG)} if it has a DML, otherwise it is called a \emph{Non-Distance Magic (NDM) graph}. In 1996, Vilfred proposed a conjecture that cylindrical grid graphs $P_m \Box C_n$ are NDM for $m \geq 2$, $n \geq 3$ and $m,n\in\mathbb{N}$. Recently, the authors could prove the conjecture for the case when $m$ is even by introducing neighbourhood chains of Type-1 (NC-T1) and Type-2 (NC-T2). In this paper, they introduce neighbourhood chains of Type-3 (NC-T3) and using them completely settle the conjecture and also identify families of NDM graphs.
	\end{abstract}

	\maketitle

		\section{Introduction}
		
			In 1987, Vilfred \cite{v87} defined \emph{Sigma labeling} and develped its theory \cite{v96}. In 2003 \cite{mrs}, the same labeling was independently defined as \emph{1-distance magic vertex labeling} and in a 2009 article \cite{sf} the term `Distance Magic Labeling' was used. Hereafter, we use the term `Distance Magic Labeling' since construction of Magic Squares motivated to define Sigma labeling of graphs (See pages 1, 97-99 in \cite{v96}.).
		
			 In 1996, Vilfred studied Distance Magic Labeling of cylindrical grid graphs $P_m \Box C_n$ for $m \geq 2$, $n \geq 3$ and $m,n\in\mathbb{N}$ and proved that $P_2 \Box C_n$ for $n \geq 3$, and $P_k \Box C_n$ for $n$ = 3,4 and $k \geq 2$ are Non-Distance Magic (NDM) and proposed the following  conjecture in \cite{v96} (pages 11-13).
	\begin{conj}  \label{c1} {\rm \cite{v96} For $m \geq 2$, $n \geq 3$ and $m,n\in\mathbb{N}$, cylindrical grid graphs $P_m \Box C_n$ are Non-Distance Magic (NDM). } 
		\end{conj}	 
		 Recently, the authors could prove the above conjecture for the case when $m$ is even by introducing neighbourhood chains of Type-1 (NC-T1) and Type-2 (NC-T2) \cite{vsj}. In this paper,  neighbourhood chains of Type-3 (NC-T3) is defined and using them, the conjecture is completely settled. We also obtain families of NDM graphs by the presence of NC-T3 in these graphs.  
		
		Through out this paper, we consider only finite undirected simple graphs and for all basic ideas in graph theory, we follow \cite{dw01}. We present here a few basic definitions and results which are needed in this paper.
		
		\begin{definition} \cite{v87} \quad	Let $G$ be a simple graph of order $n$. Then a  bijective mapping $f: V(G) \rightarrow \left\{1,2,\ldots,n\right\} $ is called a \emph{Distance Magic Labeling (DML)} or \emph{$\Sigma$-labeling (Sigma labeling) of $G$} if ~ ${\sum_{v \in N(u)}} f(v) $ = $S$ is a constant for all  $u\in V(G)$ where $N(u)$ = $\{v\in V(G): uv\in E(G)\}$. Graph $G$ is called a \emph{Distance Magic Graph (DMG)} if it has a DML, otherwise it is called a \emph{Non-Distance Magic (NDM) Graph} and we call $S$ as \emph{Distance Magic constant}.
		 \end{definition}	  
		
		\begin{definition} \cite{vsj} \quad	 Let $N = \{N(u_{i})\}^k_{i=1}$ be a sequence of neighbourhood (nbh) in a graph $G$. The \emph{nbh sequence graph of} $N$ in $G$  denoted by $NSG[N]_G$ or simply $NSG[N]$ and is defined as the union of  all  induced subgraphs of closed nbh $N[u_{i}]$, $i$ = 1 to $k$, $k \geq 2$. i.e., $NSG[N]$ = $\bigcup^k_{i=1} {<N[u_{i}]>}$ which is a subgraph of $G$, $k\in\mathbb{N}$. 
			 \end{definition}
		 
		 	\begin{definition} \cite{vsj} \quad	 A nbh sequence $\{N(u_{i})\}^k_{i=1}$ in a graph $G$ is called a  \emph{nbh walk} if consecutive terms(nbhs) are having common point(s), $u_i\in V(G)$, $1 \leq i \leq k$ and $k \geq 2$. We denote the nbh walk $\{N(u_{i})\}^k_{i=1}$ as $N_1 N_2 \cdots  N_k$ where $N_i$ = $N(u_i)$, $u_i\in V(G)$, $k \geq 2$ and $i\in\mathbb{N}$. 
		 \end{definition}
		 \begin{definition} \cite{vsj} \quad	A nbh walk $N_1 N_2 \cdots  N_k$ is called a \emph{nbh trail} if every pair of non-consecutive terms(nbhs) of the nbh sequence have at the most one common vertex, except the first and the last terms. That is $|N_i \cap N_j|$ $\leq$ 1 for every $i,j$ $\ni$ $1 < |i-j| < k-1$, $1 \leq i,j \leq k$. 
		 \end{definition}
		 
		 \begin{definition}\cite{vsj} \quad	A nbh trial $N_1 N_2 \cdots  N_k$ is called a \emph{nbh path} or a \emph{nbh chain (NC)} if non-consecutive terms(nbhs) of the nbh sequence have no common vertex, except the first and the last terms. 
		 	
		 Nbh chain $N_1 N_2 \cdots  N_k$ is said to be of \emph{length} $k$, $k \geq 2$ and $k\in\mathbb{N}$.				
		 \end{definition}				
		 \begin{definition}\cite{vsj} \quad {\rm[\textit{Neighbourhood Chain of Type 1 (NC-T1)}]}
		 	
		 Let $n \geq 2 $,  $m \geq n+1$, $\{u_{1},u_{2},\ldots,u_{n}\}$ and $\{v_{1},v_{2},\ldots,v_{m}\}$ be two disjoint subsets of the vertex set $V(G)$ of a graph $G$ $\ni$ $u_{x}u_{y}\notin E(G)$ for every $x$ and $y$, and $N_{i}$ = $N(u_i)$ = $\{v_j\in V(G): u_{i}v_j\in E(G)$, $ 1\leq j\leq m \}$ $\neq$ $\emptyset$, $1 \leq i,x,y \leq n$ and $m,n,x,y \in \mathbb{N}$. Let $NC$ = $N_1 N_2 \cdots N_{n}$ be a neighbourhood chain in $G$. Then the nbh chain $NC$ is said to be of Type 1 (NC-T1) if it satisfies the following properties:
		 \begin{enumerate}
		 	\item [\rm (i) ] $N_{1}\setminus N_{2} = \left\{v_{1}\right\}$ and $N_{n} \setminus N_{n-1}$ = $\left\{v_{n}\right\}$;
		 		
		 	\item [\rm (ii) ]  $N_{i+1} \setminus N_{i} \subset N_{i+2}$ for $i = 1,2,\ldots,n-2$, and $\vert N_{i} \setminus N_{i+1} \vert, \vert N_{i+1} \setminus N_{i} \vert$ $\geq 1$ for $i = 1,2,\ldots,n-1$. 
		 \end{enumerate}	
		 \end{definition}
		 
		 \begin{definition}\cite{vsj} \quad \label{d3} {\rm[\textit{Type 2 Neighbourhood Chains (NC-T2)}]}
		 	
		 	Let $G$ be a graph and $N_1 N_2 \cdots N_{n}$ and  $N_{n+1} N_{n+2} \cdots N_{2n}$ be two nbh chains, each of Type-1 (NC-T1) in $G$, $n \geq 2$ and $ n\in \mathbb{N}$. Then the two nbh chains of Type-1 are called a Type 2 Neighbourhood Chains (NC-T2) if they satisfy the following properties.
		 	\begin{enumerate}
		 		\item [\rm (i) ]  $N_{1} \cap N_{2} \cap N_{n+1} = N_{n+1} \setminus N_{n+2} \neq \emptyset$ and $ N_{n} \cap N_{2n} \cap N_{2n-1} = N_{n} \setminus N_{n-1} \neq \emptyset $;
		 		\item [\rm (ii) ]  $ N_i \cap N_{i+1} \cap N_{n+i} \cap N_{n+i-1} \neq \emptyset$ for $i = 2,3,\ldots,n-1$.
		 	\end{enumerate}
		 \end{definition}	

	 	\begin{definition}\cite{dw01} The \emph{cartesian product of graphs} $G$ and $H$, written $G \Box H$, is the
	 	graph with vertex set $V(G)$ $\times$ $V(H)$ specified by putting $(u, v)$ adjacent to
	 	$(u', v')$ if and only if (i) $u = u'$ and $vv'\in E(H)$ or (ii) $v = v'$ and $uu'\in E(G)$.
	 \end{definition}		 	
		
	\begin{definition}\cite{dw01} The \emph{union of two disjoint graphs} $G$ and $H$, written $G \cup H$, is
		the graph with vertex set $ V(G) \cup V(H)$ and edge set $E(G) \cup E(H)$.
	\end{definition}
	
	\begin{theorem}  \label{a0} {\rm {[Theorem 2.6 in \cite{v96}]}\quad For $m \geq 2$, graphs $ P_m \Box C_3 $  and  $ P_m \Box C_4 $ are NDM, $ m \in \mathbb{N}$. \hfill $\Box$}
\end{theorem}

\begin{theorem}\cite{vsj} \quad \label{t1} {\rm For $m \geq 2$, $n \geq 3$, graph $P_{m} \Box C_{n}$ contains a Type-2 nbh chains, $m,n \in \mathbb{N}$.\hfill $\Box$}
\end{theorem}

\begin{theorem}\cite{vsj}\quad \label{t2} {\rm Let $G$ be a graph containing Type-1 nbh chain of length $2n$, $n\in\mathbb{N}$. Then $G$ is NDM.	\hfill $\Box$}
\end{theorem} 	

\begin{cor}\cite{vsj} \quad \label{c18}  {\rm For $n \geq 3$ and $n,k \in \mathbb{N}$, graphs $P_{2k} \Box C_{n}$ are NDM.\hfill $\Box$}	
\end{cor}
		
In this paper, we define neighbourhood chains of Type-3 (NC-T3) and using them, Conjecture \ref{c1} is completely settled. We also obtain families of NDM graphs by the presence of NC-T3 in these graphs.
	
Effort to settle Conjecture \ref{c1} is the motivation for this research work.
	
	\section{On Conditions To Identify Non-Distance Magic (NDM) Graphs} 
	
	While trying to prove the conjecture on NDM of cylindrical grid graphs, we observe that we can identify a given graph $G$ as a Non-Distance Magic graph if $G$ has four distinct vertices with three conditions on their neighbourhoods as given in the following theorem. This theorem is used to completly settle the conjecture. 

	\begin{theorem}\quad \label{t2.1} {\rm 
		Let $ G $ be a simple graph and $u_1,u_2,w_1$  and $w_2$ be distinct vertices of $G$ such that  $ N(u_1) \cap N(u_2),  N(w_1) \cap N(w_2)$ $\neq$ $\emptyset$. Let the following three conditions be satisfied on these nbhs. Then $G$ is NDM.
			\begin{enumerate}
			\item [\rm (i) ] $ N(u_1)\setminus N(u_2) $, $ N(u_2)\setminus N(u_1) $, $ N(w_1)\setminus N(w_2) $,
			$ N(w_2)\setminus N(w_1) $ $ \neq \emptyset $;
			\item [\rm (ii) ]   $ N(w_1)\setminus N(w_2) \subset N(u_1)\setminus N(u_2)$ and $ N(u_2)\setminus N(u_1) \subset N(w_2)\setminus N(w_1)$ 
			 and
			\item [\rm (iii) ] $\vert (N(u_1)\setminus N(u_2))\setminus (N(w_1)\setminus N(w_2))\vert$ 
			
			\hfill + $\vert(N(w_2)\setminus N(w_1))\setminus (N(u_2)\setminus N(u_1))\vert $ $ \geq 1 $.
		\end{enumerate}
	}
	\end{theorem}
	\begin{proof}
		Suppose $G$ is 
		DM. Let $f$ be a DML of $G$ with DM constant $S$. By the definition of DML, we get,
		\begin{align}
			S~=\sum_{v_j \in  N(u_1)}f(v_j)=~\sum_{v_j \in N(u_2) }f(v_j)~=~\sum_{v_j \in  N(w_1)}f(v_j)~=~\sum_{v_j \in N(w_2)}f(v_j). \notag
		\end{align} 
	Since $ N(u_1) \cap N(u_2) \neq \emptyset $, 
		\begin{align}
		\sum_{v_j \in N(u_1)}f(v_j)=\sum_{v_j \in N(u_2) }f(v_j)~~~ \implies\notag
			\end{align} 
\begin{align}
\sum_{v_j \in N(u_1)\setminus N(u_2)}f(v_j)=\sum_{v_j \in N(u_2)\setminus N(u_1)}f(v_j). \label{a} 
	\end{align}
	Also, $ N(w_1)\setminus N(w_2) \subset N(u_1)\setminus N(u_2)$, (\ref{a}) becomes,
\begin{align}
\sum_{v_j \in ((N(u_1)\setminus N(u_2))\setminus (N(w_1)\setminus N(w_2)))}f(v_j) +\sum_{v_j \in N(w_1)\setminus N(w_2)}f(v_j)=\notag\\\sum_{v_j \in N(u_2)\setminus N(u_1)}f(v_j).\label{c} 
	\end{align} 
	Similarly,  $ N(w_1) \cap N(w_2) \neq \emptyset $,
		\begin{align}
		\sum_{v_j \in N(w_1)}f(v_j)=\sum_{v_j \in N(w_2) }f(v_j)\Rightarrow~~~ ~~~~~~~~~~\notag\\ \sum_{v_j \in N(w_1)\setminus N(w_2)}f(v_j)=\sum_{v_j \in N(w_2)\setminus N(w_1)}f(v_j) \label{b}. 
	\end{align} 
	Since $N(u_2)\setminus N(u_1) \subset N(w_2)\setminus N(w_1)$, (\ref{b}) becomes,
	\begin{align}
		\sum_{v_j \in N(w_1)\setminus N(w_2)}f(v_j)=\sum_{v_j \in ((N(w_2)\setminus N(w_1))\setminus (N(u_2)\setminus N(u_1)))}f(v_j) +\notag\\\sum_{v_j \in N(u_2)\setminus N(u_1)}f(v_j).\label{d}
	\end{align} 
	 From (\ref{c}) \& (\ref{d}), we get, \small
	 \begin{align}
	 	\sum_{v_j \in ((N(u_1)\setminus N(u_2))\setminus (N(w_1)\setminus N(w_2)))}f(v_j) +\sum_{v_j \in ((N(w_2)\setminus N(w_1))\setminus (N(u_2)\setminus N(u_1)))}f(v_j) = 0, \notag 
 	\end{align} 
	which is a contradiction to condition (iii), since $ ((N(u_1)\setminus N(u_2))\setminus (N(w_1)\setminus N(w_2))) \cup ((N(w_2)\setminus N(w_1))\setminus (N(u_2)\setminus N(u_1)))$  contains atleast one element whose labeling is $ \geq 1$.\\
	 Hence the result. 
	 	\end{proof}

We illustrate above theorem with the following two examples.
\begin{exm} \quad Let $G $ be a graph as given in Figure 1. Then $ G $ is NDM.
	{\rm \begin{center}
		\begin{tikzpicture}  
			[scale=.4,auto=center,every node/.style={draw,circle}] \node (1) at (-8,2)[circle,fill=blue!30]  {\tiny $u_{1}$};
			
			\node (2) at (-2,2)[circle,fill=blue!30] {\tiny$u_{2}$};
			\node (3) at (-11,4)  {\tiny $v_{1}$};
			\node (c) at (-11,0)  {\tiny$v_{2}$};
			\node (4) at (-5,0) {\tiny$v_{3}$};
			\node (11) at (1,0){\tiny$v_{4}$};
			\node (13) at (-8,-2) [circle,fill=blue!30]  {\tiny$w_{1}$};
			\node (14) at (-2,-2) [circle,fill=blue!30]  {\tiny$w_{2}$};
			\node (16) at (-5,-4) {\tiny$v_{5}$};
			\node (18) at (1,-4) {\tiny$v_{6}$};	
			
			\draw (1) -- (3);
			\draw (1) -- (c);
			\draw (13) -- (c);
			
			\draw (1) -- (4);
			\draw (2) -- (4);
			\draw (2) -- (11);
			\draw (11) -- (14);
			\draw (13) -- (16);
			\draw (14) -- (18);	
			\draw (14) -- (16);
		\end{tikzpicture}\\  
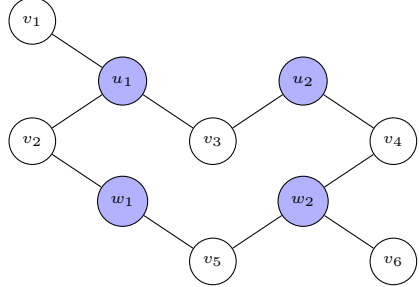
		Figure 1.\label{fig 1}  
	\end{center}

For more clarity, we use Venn diagram to represent the neighbourhoods of the vertices $u_1,$ $u_2,$ $w_1,$ $w_2 $ of the graph given in Figure 1 and present it in Figure 2. Here, the neighbourhoods of the vertices $ u_1,$ $u_2,$ $w_1,$ $w_2 $ of the graph are represented with colors blue, green, brown, red, respectively. And these neighbourhoods satisfy the three conditions given in Theorem \ref{t2.1}. See Figure 2. Thereby the graph $G$ is NDM.}	
	\begin{center}
		\begin{tikzpicture}[thick,
			set/.style = { circle, minimum size = .061cm}]
			
			\draw[green] (3,0.5) ellipse (1.9cm and 0.5cm);
				\draw[brown] (3,-0.8) ellipse (1.7cm and 0.5cm);
			
			\draw[blue] (1,0) ellipse (2.3cm and 1.4cm);
				\begin{scope}[red]
				\clip(2.15,-2)rectangle(7,2);
				\draw (1.4,0) circle [x radius=1.9cm, y radius=1.37cm];
			\end{scope}
			
			\begin{scope}[red]
			\clip(2.15,-2)rectangle(7,2);
			\draw (1.4,0) circle [x radius=1.92cm, y radius=1.36cm];
		\end{scope}
			\begin{scope}[blue]
				\clip(2.15,-2)rectangle(7,2);
				\draw (1.4,0) circle [x radius=1.86cm, y radius=1.35cm];
			\end{scope}

			\begin{scope}[red]
				\clip(2.1,-1.7)rectangle(7,2);
				\draw (4,0) circle [x radius=3cm, y radius=15.8mm];
			\end{scope}
			\node[black] at (0,0) {$v_1$};
				\draw[->] [blue](-.6,-.5) to [out=50,in=55] (-.3,-2);
			\node[blue] at (-.5,-2.4) {$N(u_1)$};
			
				\draw[->] [green](3.6,0.6) to [out=50,in=55] (3,1.8);
			\node[green] at (2.5,1.9) {$N(u_2)$};
				\draw[->] [brown](3.3,-0.9) to [out=50,in=55] (3.6,-2);
			\node[brown] at (3.4,-2.4) {$N(w_1)$};
		
			\draw[->] [red](6.7,0) to [out=35,in=120] (7.6,-0.3);
			\node[red] at (7.6,-0.6) {$N(w_2)$};
				\node[black] at (2,0.6) {$v_3$};
					\node[black] at (2,-.9) {$v_2$};
				\node[black] at (4,0.6) {$v_4$};
				\node[black] at (4,-0.9) {$v_5$};
				\node[black] at (6,0) {$v_6$};
					\draw [olive] (-2.5,-3) rectangle (8.5,2.5);
		\end{tikzpicture}\\
		\vspace{0.2cm} Figure $2$. Venn diagram of Figure $ 1 $. \label{fig 2}			
	\end{center}
\end{exm}
\begin{exm} \quad {\rm Let $G $ be a graph as given in Figure 3. Then $ G $ is NDMG.

		\begin{center}
	\begin{tikzpicture}  
		[scale=.58,auto=center,every node/.style={draw,circle}] 
		
		\node (1) at (-8,2)[circle,fill=blue!30]  {$u_{1}$};	
		\node (2) at (-2,2)[circle,fill=blue!30] {$u_{2}$};
		\node (3) at (-10,4)  {$v_{1}$};
		\node (a) at (-10,2)  {$v_{2}$};
		\node (b) at (-10,0)  {$v_{3}$};
		\node (c) at (-7,0)  {$v_{4}$};
		\node (4) at (-5,1) {$v_{5}$};
		
		\node (10) at (-5,-1) {$v_{6}$};
		\node (11) at (-3.5,0){$v_{7}$};
		\node (12) at (-5,4){$v_{9}$};
		\node (d) at (-0.5,0){$v_{8}$};
		\node (13) at (-8,-2) [circle,fill=blue!30]  {$w_{1}$};
		\node (14) at (-2,-2) [circle,fill=blue!30]  {$w_{2}$};
		\node (16) at (-4.8,-4) {$v_{10}$};
		\node (17) at (-2,-4) {$v_{11}$};
		\node (18) at (0,-3) {$v_{12}$};	
		
		\draw (1) -- (3);
		\draw (1) -- (a);
		\draw (1) -- (b);
		\draw (1) -- (c);
		\draw (13) -- (b);
		\draw (13) -- (c);
		\draw (1) -- (12);
		\draw (2) -- (12);
		
		\draw (1) -- (4);
		\draw (13) -- (10);
		\draw (2) -- (4);
		\draw (14) -- (10);
		\draw (2) -- (11);
		\draw (2) -- (d);
		\draw (14) -- (d);
		
		\draw (11) -- (14);
		\draw (13) -- (16);
		\draw (14) -- (17);
		\draw (14) -- (18);	
		\draw (14) -- (16);
	\end{tikzpicture}\\  
	Figure $ 3 $.  
\end{center}
\begin{center}
	\begin{tikzpicture}[thick,
		set/.style = { circle, minimum size = .061cm}]
		
		\draw[green] (3,0.5) ellipse (1.9cm and 0.5cm);
		\draw[brown] (3,-0.8) ellipse (1.7cm and 0.5cm);
		
		\draw[blue] (1,0) ellipse (2.3cm and 1.4cm);
		\begin{scope}[red]
			\clip(2.15,-2)rectangle(7,2);
			\draw (1.4,0) circle [x radius=1.9cm, y radius=1.37cm];
		\end{scope}
		
		\begin{scope}[red]
			\clip(2.15,-2)rectangle(7,2);
			\draw (1.4,0) circle [x radius=1.92cm, y radius=1.36cm];
		\end{scope}
		\begin{scope}[blue]
			\clip(2.15,-2)rectangle(7,2);
			\draw (1.4,0) circle [x radius=1.86cm, y radius=1.35cm];
		\end{scope}

		\begin{scope}[red]
			\clip(2.1,-1.7)rectangle(7,2);
			\draw (4,0) circle [x radius=3cm, y radius=15.8mm];
		\end{scope}
		\node[black] at (0,0.6) {$v_1$};
		\node[black] at (0.5,-0.1) {$~v_2$};
		\draw[->] [blue](0,-0.9) to [out=50,in=55] (0,-2);
	\node[blue] at (.5,-2.4) {$N(u_1)$};
		\draw[->] [green](2.3,0.6) to [out=50,in=55] (2,1.8);
	\node[green] at (1.3,1.9) {$N(u_2)$};
	\draw[->] [brown](3.3,-0.9) to [out=50,in=55] (3.6,-2);
	\node[brown] at (3.4,-2.4) {$N(w_1)$};
	
	\draw[->] [red](6.7,0) to [out=35,in=120] (7.6,-0.3);
	\node[red] at (7.6,-0.6) {$N(w_2)$};
\node[black] at (2,0.6) {$v_5 $};
\node[black] at (2.5,0.3) {$~v_9$};
		\node[black] at (2,-.9) {$v_3~$};
		\node[black] at (2.5,-.6) {$~v_4$};
		\node[black] at (3.5,0.7) {$v_7~$};
			\node[black] at (4,0.4) {$~v_8$};
		\node[black] at (3.5,-0.6) {$v_6~$};
		\node[black] at (4,-0.9) {$~v_{10}$};
		\node[black] at (5.5,0.9) {$v_{11}~$};
		\node[black] at (5.8,-0.5) {$~v_{12}$};
			\draw [olive] (-2.5,-3) rectangle (8.5,2.5);
	\end{tikzpicture}\\
	\vspace{0.2cm} Figure $4$. Venn Diagram of Figure $3 $. 			
\end{center}

Similar to Figure 2, the Venn diagram corresponding to Figure 3 is given in Figure 4. And these neighbourhoods satisfy all the conditions given in Theorem \ref{t2.1}. See Figure 4. Thereby the graph $G$ is NDM.}
\end{exm}

	\section{Type-3 nbh chains and existence of NDM graphs}
	
In this section, we introduce a nbh chains of Type 3 (NC-T3), prove that an NC-T2 with odd length is an NC-T3 and completly settle Conjecture \ref{c1}.

\begin{definition} \quad \label{d3.1} {\rm[\textit{Type 3 Neighbourhood Chains (NC-T3)}]}
		
		Let $G$ be a graph, $n_1,n_2 \geq 2 $, $n_3 \geq \max \{n_1,n_2\} + 3$ and $n_1,n_2,n_3\in\mathbb{N}$. For $k$ = 1,2,3, let $V_k$ = $\{u^k_{1},u^k_{2},\ldots,u^k_{n_k}\}$ $\subset$ $V(G)$, $V_1$ and $V_2$ be stable sets in $G$, and $V_i$ $\cap$ $V_j$ $=$ $\emptyset$ for $i \neq j$, $ 1\leq i,j \leq 3$. For $j$ = 1,2, let $N^j_{i}$ = $N(u^j_i)$ = $\{u^3_k\in V(G): u^j_{i}u^3_k\in E(G)$, $ 1\leq k\leq n_3 \}$ $\neq$ $\emptyset$, $1 \leq i \leq n_j$. For $i$ = 1,2, let $NC_i$ = $N^i_1 N^i_2 \cdots N^i_{n_i}$ be a neighbourhood chain in $G$. Then, $NC_1$ and $NC_2$ are said to be \emph{nbh chains of Type 3} $(NC-T3)$ if they satisfy the following conditions:
		
		
		\begin{enumerate}
			\item [\rm (i) ] $N^j_{i+1} \setminus N^j_{i} \subset N^j_{i+2}$ for $i$ = 1 to $n_j-2$ and $j$ = 1,2 and
			\\
			$\vert N^j_{i} \setminus N^j_{i+1} \vert, \vert N^j_{i+1} \setminus N^j_{i} \vert$ $\geq 1$ for $i$ = 1 to $n_j-1$ and $j$ = 1,2;
			
			\item [\rm (ii) ] Let $q = \min \{n_1,n_2\}$, $1 \leq i \leq n_1-1$,  $1 \leq j \leq n_2-1$ and $0 \leq r \leq \left\lfloor{\frac{q}{2}}\right\rfloor -1 
			$. Then, for atleast one set of values of $ i, j,$ and $r $, 
			
			$ N^2_{j} \setminus N^2_{j+1} \subset N^1_{i} \setminus N^1_{i+1} $, $N^1_{i+1+2r} \setminus 
			N^1_{i+2r} \subset  N^2_{j+1+2r} \setminus N^2_{j+2r} $ and 
			
			$\vert (N^1_{i} \setminus N^1_{i+1})\setminus (N^2_{j} \setminus N^2_{j+1})\vert $ + $\vert  ( N^2_{j+1+2r} \setminus N^2_{j+2r}) \setminus (N^1_{i+1+2r} \setminus 
			N^1_{i+2r}) \vert \geq 1$.
		\end{enumerate}

			
			
			
			
			
	 
	\end{definition}

The following two examples illustrate existence of graphs $G$ with NC-T3.

\begin{exm} \quad Let $G $ be a graph as given in Figure 5. We follow the notations given in the above definition of NC-T3. Let $NC_1$ = $N^1_1N^1_2N^1_3 $ and $NC_2$ = $N^2_1N^2_2N^2_3 $ where $N^j_{i}$ = $N(u^j_i)$ for $i$ = 1,2,3 and $j$ = 1,2. Then, $NC_1$ and $NC_2$ form an NC-T3 in $G$.
	
Consider, graph $G$ as given in Figure 5. Venn diagram of  neighbourhoods $N^1_1$, $N^1_2$, $N^1_3$, $N^2_1$, $N^2_2$, $N^2_3$ of $G$ of the vertices $u^1_1$, $u^1_2$, $u^1_3$, $u^2_1$, $u^2_2$, $u^2_3$ is given in Figure $6$. See Figure 6. 
	
	{\rm 	\begin{center}
			\begin{tikzpicture}  
				[scale=.6,auto=center,every node/.style={draw,circle}] 
				
				\node (1) at (-10,2)[circle,fill=blue!30]  {$u^1_{1}$};
				\node (2) at (-5,2)[circle,fill=blue!30] {$u^1_{2}$};
				\node (3) at (-12.5,4)  {$u^3_{1}$};
				\node (a) at (-12.5,2)  {$u^3_{2}$};
				\node (b) at (-12.5,0)  {$u^3_{3}$};
				\node (c) at (-7.5,0)  {$u^3_{4}$};
				\node (4) at (-7.5,2) {$u^3_{5}$};
				\node (7) at (0,2) [circle,fill=blue!30] {$u^1_{3}$};
				
				\node (10) at (-7.5,-2) {$u^3_{6}$};
				\node (11) at (-2.5,0){$u^3_{7}$};
				\node (12) at (2.5,0){$u^3_{8}$};
				\node (d) at (2.5,4){$u^3_{9}$};
				\node (13) at (-10,-2) [circle,fill=blue!30]  {$u^2_{1}$};
				\node (14) at (-5,-2) [circle,fill=blue!30]  {$u^2_{2}$};
				\node (15) at (0,-2)  [circle,fill=blue!50] {$u^2_{3}$};
				\node (16) at (-7.5,-4) {$u^3_{10}$};
				\node (17) at (-2.5,-4) {$u^3_{11}$};
				\node (18) at (2.5,-4) {$u^3_{12}$};	
				
				\draw (1) -- (3);
				\draw (1) -- (a);
				\draw (1) -- (b);
				\draw (1) -- (c);
				\draw (13) -- (b);
				\draw (13) -- (c);
				\draw (d) -- (7);
				\draw (15) -- (11);
				
				\draw (1) -- (4);
				\draw (13) -- (10);
				\draw (2) -- (4);
				\draw (14) -- (10);
				\draw (2) -- (11);
				\draw (7) -- (11);
				
				\draw (7) -- (12);
				\draw (11) -- (14);
				\draw (13) -- (16);
				\draw (14) -- (17);
				\draw (15) -- (18);	
				\draw (14) -- (16);
				\draw (15) -- (17);
			\end{tikzpicture}\\
			
			\vspace{.1cm}  
			Figure  $5$. \hspace{2cm}~ 
				\end{center}
			
			\begin{center}
		\begin{tikzpicture}[ thick, main/.style = {draw}]
			\begin{scope}[magenta]
				\clip(0,0)rectangle(5,3);
				\draw (3.5,1.5) circle [x radius=3cm, y radius=15mm];
			\end{scope}
			\draw[magenta] (4.97,0.2) -- (4.97,2.8);
			\begin{scope}[orange]
				\clip(0,0)rectangle(7.47,3);
				\draw (6.,1.5) circle [x radius=3.3cm, y radius=13.65mm];
			\end{scope}
			\draw[orange] (7.47,0.27) -- (7.47,2.73);
			\begin{scope}[black]
				\clip(5.07,1)rectangle(7.43,3);
				\draw (6,1.43) circle [x radius=3.3cm, y radius=13.65mm];
			\end{scope}
			\begin{scope}[red]
				\clip(5.18,1)rectangle(7.39,3);
				\draw (6,1.38) circle [x radius=3.3cm, y radius=13.65mm];
			\end{scope}
			\begin{scope}[cyan]
				\clip(5.04,0)rectangle(10,3);
				\draw (6,1.5) circle [x radius=4cm, y radius=14mm];
			\end{scope}
			\draw[cyan] (5.04,0.15) -- (5.04,2.85);
			


			\draw[green] (2,-0.5) ellipse (1cm and 2cm);
			\begin{scope}[black]
				\clip(0.8,0.55)rectangle(2,-4);
				\draw (1.95,-0.5) ellipse (1cm and 2cm);
			\end{scope}
			
			\begin{scope}[black]
				\clip(1.1,-2)rectangle(3,1.4);
				\draw (3.45,1.5) circle [x radius=3.05cm, y radius=15.4mm];
			\end{scope}	
			\begin{scope}[black]
				\clip(2.9,-2)rectangle(5.1,2);
				\draw (5,1) circle [x radius=3cm, y radius=14mm];
			\end{scope}			
			
			\draw[black] (2,-2.5) -- (6.5,-2.5);
			\draw[black] (7.4,-0.5) -- (7.4,2.7);
			\draw[black] (5.09,-0.4) -- (5.09,2.75);
			
			\begin{scope}[black]
				\clip(6.4,-0.5)rectangle(8,-4);
				\draw (6.4,-0.5) ellipse (1cm and 2cm);
			\end{scope}
			\begin{scope}[red]
				\clip(2,-0.1)rectangle(4,-4);
				\draw (2.05,-0.5) ellipse (1cm and 2.05cm);
			\end{scope}
			\begin{scope}[red]
				\clip(3,-2)rectangle(5.18,1);
				\draw (5.15,0.9) circle [x radius=3cm, y radius=14mm];
			\end{scope}
			\draw[red] (2,-2.56) -- (7.35,-2.56);
			\draw[red] (7.33,-0.47) -- (7.33,2.65);
			\draw[red] (5.16,-0.5) -- (5.16,2.73);
			\begin{scope}[red]
				\clip(7.3,0)rectangle(10,-4);
				\draw (5.9,-1.5) circle [x radius=4cm, y radius=11.5mm];
			\end{scope}
			\node at (1.5,2) {$u_1^3$};
			\node at (2.3,2.5) {$u_2^3$};
			\node at (3.9,1.6) {$u_5^3$};
			\node at (1.6,1) {$u_3^3$};
			\node at (2.3,0.5) {$u_4^3$};
			\node at (5.9,1.6) {$u_7^3$};
			\node at (8,1) {$u_8^3$};
			\node at (9,2) {$u_9^3$};
			\node at (1.4,-1) {$u_6^3$};
			\node at (2.3,-0.5) {$u_{10}^3$};
			\node at (5.5,-1.2) {$u_{11}^3$};
			\node at (8.5,-1.1) {$u_{12}^3$};
			\draw [olive] (-0.5,-3.8) rectangle (11,3.8);
			\draw[->] [magenta](1,2) to [out=35,in=120] (0,1.5); 
			\node [magenta]  at (0,1.3) {$N_{1}^1$};
			
			\draw[->] [orange](4.2,2) to [out=50,in=55] (4,3.5);
			\node [orange]  at (3.6,3.4) {$N_{2}^1$};
			
			\draw[->] [cyan](9.2,1.5) to [out=35,in=120] (10.5,1.3); 
			\node [cyan]  at (10.5,1) {$N_{3}^1$};
			
			\draw[->] [red](9.2,-1.1) to [out=35,in=120] (10.5,-1.4); 
			\node [red]  at (10.5,-1.7) {$N_{3}^2$};
			\draw[->] [green](1.4,-1.5) to [out=35,in=120] (0.5,-1.8); 
			\node [green]  at (0.5,-2.2) {$N_{1}^2$};
			\draw[->] [black](4.2,-2.2) to [out=50,in=55] (4.5,-3);
			\node [black]  at (4.7,-3.3) {$N_{2}^2$};

		\end{tikzpicture}\\
		\vspace{.1cm}  
	Figure  $6$. Venn Diagram of Figure $5$.  
\end{center}	
	
		In Figure 6, $ N^1_1N^1_2N^1_3 $ and $ N^2_1N^2_2N^2_3 $ are two nbh chains in $G$. Also,
		\begin{align}
			N^1_2 \setminus N^1_1 &= \{u_{7}^3\} \subset \{u_{7}^3,u_{8}^3,u_{9}^3\} = N^1_3, \label{22}\\
			N^2_2 \setminus N^2_1 &= \{u_{7}^3,u_{11}^3\} \subset \{u_{7}^3,u_{11}^3,u_{12}^3\} = N^2_3, \label{23}\\
			N^2_1 \setminus N^2_2 &= \{u_{3}^3,u_{4}^3\} \subset \{u_{1}^3,u_{2}^3,u_{3}^3,u_{4}^3\} = N^1_1 \setminus N^1_2, \label{24}\\
			N^1_2 \setminus N^1_1 &= \{u_{7}^3\} \subset  \{u_{7}^3,u_{11}^3\} = N^2_2 \setminus N^2_1.\label{25}
		\end{align}
		It is clear that the nbh chains $ N^1_1N^1_2N^1_3 $ = $NC_1$ and $ N^2_1N^2_2N^2_3 $ = $NC_2$ satisfy property (i) of NC-T3 by (\ref{22}) and (\ref{23}), and property (ii) of NC-T3 by (\ref{24}) and (\ref{25}) where $i$ = 1 = $j$ and $ r$ = 0. Thus  $ NC_1 $ and $ NC_2 $ form an NC-T3 in $G$. See Figures 5 and 6.}
	\end{exm}	
	
	\begin{exm} \quad Let $ G $ be a graph as given in Figure 7. Let $NC_1$  $=  N^1_1N^1_2N^1_3N^1_4N^1_5\\N^1_6N^1_7 $ and $NC_2$ $= N^2_1N^2_2N^2_3N^2_4 $ form NC-T3 in $G$, where $N^1_{i}$ = $N(u^1_i)$ and  $N^2_{j}$ = $N(u^2_j)$  for $i
		 = 1,2,\ldots,7$ and $j$ = 1,2,3,4. {\rm
			
		\begin{center}
			\begin{tikzpicture}  
				[scale=.7,auto=center,every node/.style={draw,circle}] 
				
				\node (1) at (-10,4)[circle,fill=blue!30]  {$u_{1}^1$};
				\node (2) at (-8,4)[circle,fill=blue!30]  {$u_{2}^1$};
				\node (3) at (-6,4)[circle,fill=blue!30]  {$u_{3}^1$};
				\node (4) at (-4,4)[circle,fill=blue!30]  {$u_{4}^1$};
				\node (5)at(-2,4)[circle,fill=blue!30]{$u_{5}^1$};		
				\node (6) at (0,4)[circle,fill=blue!30]  {$u_{6}^1$};
				\node (7) at (2,4)[circle,fill=blue!30]  {$u_{7}^1$};
				
				\node (13) at (-6,-0.5)[circle,fill=blue!30]  {$u_{1}^2$};
				\node (14) at (-4,-0.5)[circle,fill=blue!30]  {$u_{2}^2$};
				\node (15) at (-2,-0.5)[circle,fill=blue!30]  {$u_{3}^2$};
				\node (16) at (0,-0.5)[circle,fill=blue!30]  {$u_{4}^2$};
				\node (a) at (-11,5.5)  {$u_{1}^3$};
				\node (b) at (-11,3)  {$u_{2}^3$};
				\node (c) at (-8.8,5.5) {$u_{3}^3$};
				\node (e) at (-6.9,5.5)  {$u_{4}^3$};
				\node (f) at (-7,2)  {$u_{5}^3$};
				\node (g) at (-5,5.5)  {$u_{6}^3$};
				\node (h) at (-5,2.5)  {$u_{7}^3$};	
				\node (i) at (-2.8,5.5) {$u_{8}^3$};
				\node (j) at (-3,2.5)  {$u_{9}^3$};
				\node (k) at (-0.8,5.5)  {$u_{10}^3$};
				\node (m) at (0.8,5.5)  {$u_{11}^3$};
				\node (n) at (1,2)  {$u_{12}^3$};
				\node (o) at (2.8,5.5)  {$u_{13}^3$};
				\node (p) at (3,2.5)  {$u_{14}^3$};	
				
				\node (g1) at (-5,1) {$u_{15}^3$};
				\node (h1) at (-5,-2.2) {$u_{16}^3$};
				\node (i1) at (-3,1)  {$u_{17}^3$};
				\node (j1) at (-3,-2.2)  {$u_{18}^3$};
				\node (k1) at (-1,1) {$u_{19}^3$};
				\node (l1) at (-1,-2.2) {$u_{20}^3$};	
				\node (n1) at (0.8,-2.2)  {$u_{21}^3$};
				
				\draw (1) -- (a);
				\draw (1) -- (b);	
				\draw (1) -- (c);
				\draw (2) -- (c);
				\draw (2) -- (e);
				\draw (3) -- (e);
				\draw (2) -- (f);
				\draw (3) -- (f);
				\draw (3) -- (g);
				\draw (4) -- (g);
				\draw (3) -- (h);
				\draw (4) -- (h);
				\draw (4) -- (i);
				\draw (5) -- (i);
				\draw (4) -- (j);
				\draw (5) -- (j);
				\draw (5) -- (k);
				\draw (6) -- (k);
				
				\draw (7) -- (m);
				\draw (6) -- (n);
				\draw (7) -- (n);
				\draw (7) -- (o);
				\draw (7) -- (p);
				\draw (13) -- (f);
				\draw (13) -- (g1);
				\draw (13) -- (h1);
				
				\draw (14) -- (g1);
				\draw (14) -- (h1);
				\draw (14) -- (i1);
				\draw (14) -- (j1);
				\draw (15) -- (i1);
				\draw (15) -- (j1);
				\draw (15) -- (k1);
				\draw (15) -- (l1);
				\draw (16) -- (k1);
				\draw (16) -- (l1);
				\draw (16) -- (n);
				\draw (16) -- (n1);
				
			\end{tikzpicture}\\  
			\vspace{0.2cm} Figure $ 7$.  
			
		\end{center}
		Consider, graph $G$ as given in Figure 7. Venn Diagram of neighbourhoods  $N^1_1$, $N^1_2$, $N^1_3$, $N^1_4$, $N^1_5$, $N^1_6$, $N^1_7$, $N^2_1$, $N^2_2$, $N^2_3$, $N^2_4$ of $G$ of the vertices $u^1_1$, $u^1_2$, $u^1_3$, $u^1_4$, $u^1_5$, $u^1_6$, $u^1_7$, $u^2_1$, $u^2_2$, $u^2_3$, $u^2_4$ is given in Figure $8$. See Figures 7 and 8.  	
	\begin{center}
		\begin{tikzpicture}[thick,
			set/.style = { circle, minimum size = .061cm}]
			\draw[yellow] (3,0) ellipse (2cm and 1cm);
			\draw[yellow] (3.01,0.01) ellipse (2cm and 1cm);
			\draw[yellow] (2.98,-0.03) ellipse (2cm and 1cm);
			\draw[black] (1,0) ellipse (1.94cm and 0.99cm);
			\begin{scope}[red]
				\clip(2,-1)rectangle(3,1.2);
				\draw (1,0) circle [x radius=2cm, y radius=10mm];
			\end{scope}
			\begin{scope}[red]
				\clip(2,-1)rectangle(3,1.2);
				\draw (1,0) circle [x radius=1.975cm, y radius=10mm];
			\end{scope}
			\begin{scope}[black]
				\clip(2,-1)rectangle(3,1.2);
				\draw (1,0) circle [x radius=1.935cm, y radius=9.8mm];
			\end{scope}
			
			\begin{scope}[red]
				\clip(2,-1.5)rectangle(3.5,1.5);
				\draw (3.1,0) ellipse (2.05cm and 1.04cm);
			\end{scope}
			\begin{scope}[red]
				\clip(2,-1.5)rectangle(3.5,1.5);
				\draw (3.1,0) ellipse (2.05cm and 1.05cm);
			\end{scope}
			\begin{scope}[red]
				\clip(3.5,-1.5)rectangle(8,1.5);
				\draw (3.9,0) ellipse (2.5cm and 1.03cm);
			\end{scope}
			\begin{scope}[red]
				\clip(3.5,-1.5)rectangle(8,1.5);
				\draw (3.9,0) ellipse (2.5cm and 1.01cm);
			\end{scope}
			\begin{scope}[red]
				\clip(3.5,-1.5)rectangle(8,1.5);
				\draw (3.9,0) ellipse (2.5cm and 1.04cm);
			\end{scope}
			\begin{scope}[cyan]
				\clip(3.5,-1)rectangle(7,1.2);
				\draw (3.1,0) circle [x radius=1.95cm, y radius=10.5mm];
			\end{scope}
			\begin{scope}[cyan]
				\clip(3.5,-1)rectangle(7,1.2);
				\draw (3.1,0) circle [x radius=1.95cm, y radius=10.7mm];
			\end{scope}
			\begin{scope}[cyan]
				\clip(3.5,-1)rectangle(7,1.2);
				\draw (3.1,0) circle [x radius=1.95cm, y radius=10.75mm];
			\end{scope}
			\begin{scope}[cyan]
				\clip(3.6,-1.5)rectangle(4.5,1.5);
				\draw (4.1,0) ellipse (2.05cm and 1.07cm);
			\end{scope}
			
			\begin{scope}[cyan]
				\clip(3.6,-1.5)rectangle(4.5,1.5);
				\draw (4.1,0) ellipse (2.05cm and 1.08cm);
			\end{scope}
			\begin{scope}[cyan]
				\clip(3.6,-1.5)rectangle(4.5,1.5);
				\draw (4.1,0) ellipse (2.05cm and 1.085cm);
			\end{scope}
			\begin{scope}[cyan]
				\clip(4.5,-1.5)rectangle(9,1.5);
				\draw (4.8,0) ellipse (2.9
				cm and 1.085cm);
			\end{scope}
			\begin{scope}[cyan]
				\clip(4.5,-1.5)rectangle(9,1.5);
				\draw (4.8,0) ellipse (2.83cm and 1.08cm);
			\end{scope}
			\begin{scope}[cyan]
				\clip(4.5,-1.5)rectangle(9,1.5);
				\draw (4.8,0) ellipse (2.82cm and 1.088cm);
			\end{scope}
			\begin{scope}[cyan]
				\clip(4.5,-1.5)rectangle(9,1.5);
				\draw (4.8,0) ellipse (2.81cm and 1.089cm);
			\end{scope}
			\begin{scope}[cyan]
				\clip(4.5,-1.5)rectangle(9,1.5);
				\draw (4.8,0) ellipse (2.8cm and 1.09cm);
			\end{scope}
			\begin{scope}[green]
				\clip(4.5,-1.5)rectangle(8,1.5);
				\draw (3.9,0) ellipse (2.55cm and 1.06cm);
			\end{scope}
			\begin{scope}[green]
				\clip(4.5,-1.5)rectangle(8,1.5);
				\draw (3.9,0) ellipse (2.57cm and 1.08cm);
			\end{scope}
				
			\begin{scope}[green]
				\clip(4.5,-1.5)rectangle(5.9,1.5);
				\draw (5.3,0) ellipse (2.05cm and 1.09cm);
			\end{scope}
			\begin{scope}[green]
				\clip(4.5,-1.5)rectangle(5.9,1.5);
				\draw (5.3,0) ellipse (2.05cm and 1.1cm);
			\end{scope}
			\begin{scope}[green]
				\clip(5.7,-1.5)rectangle(10,1.5);
				\draw (6.2,0) ellipse (2.9cm and 1.05cm);
			\end{scope}
			\begin{scope}[green]
				\clip(5.7,-1.5)rectangle(10,1.5);
				\draw (6.2,0) ellipse (2.9cm and 1.03cm);
			\end{scope}
			\begin{scope}[green]
				\clip(5.7,-1.5)rectangle(10,1.5);
				\draw (6.2,0) ellipse (2.9cm and 1.07cm);
			\end{scope}
			\begin{scope}[green]
				\clip(5.7,-1.5)rectangle(10,1.5);
				\draw (6.2,0) ellipse (2.9cm and 1.08cm);
			\end{scope}
			\begin{scope}[orange]
				\clip(5.7,-1.5)rectangle(9,1.5);
				\draw (4.7,0) ellipse (3cm and 1.15cm);
			\end{scope}
			\begin{scope}[orange]
				\clip(5.7,-1.5)rectangle(9,1.5);
				\draw (4.7,0) ellipse (3cm and 1.2cm);
			\end{scope}
			\begin{scope}[orange]
				\clip(5.7,-1.5)rectangle(9,1.5);
				\draw (4.68,0) ellipse (3cm and 1.2cm);
			\end{scope}
			\begin{scope}[orange]
				\clip(5.7,-1.5)rectangle(7,1.5);
				\draw (6.3,0) ellipse (3.5cm and 1.13cm);
			\end{scope}
			\begin{scope}[orange]
				\clip(5.7,-1.5)rectangle(7,1.5);
				\draw (6.3,0) ellipse (2.9cm and 1.14cm);
			\end{scope}
			\begin{scope}[orange]
				\clip(5.7,-1.5)rectangle(7,1.5);
				\draw (6.3,0) ellipse (2.9cm and 1.119cm);
			\end{scope}
			\begin{scope}[orange]
				\clip(6.85,-1.5)rectangle(11,1.5);
				\draw (7.54,0) ellipse (2.9cm and 1.11cm);
			\end{scope}
			\begin{scope}[orange]
				\clip(6.9,-1.5)rectangle(11,1.5);
				\draw (7.55,0) ellipse (2.9cm and 1.12cm);
			\end{scope}
			\begin{scope}[pink]
				\clip(7,-1.5)rectangle(10,1);
				\draw (6.27,0.03) ellipse (2.9cm and 1.135cm);
			\end{scope}
			\begin{scope}[pink]
				\clip(7,-1.5)rectangle(10,1.3);
				\draw (6.27,0) ellipse (2.89cm and 1.13cm);
			\end{scope}
			\begin{scope}[pink]
				\clip(7,-1.5)rectangle(10,1.5);
				\draw (6.27,0) ellipse (2.89cm and 1.15cm);
			\end{scope}
			\begin{scope}[pink]
				\clip(7,-1.5)rectangle(8.5,1.5);
				\draw (7.8,0) ellipse (2.9cm and 1.145cm);
			\end{scope}
			\begin{scope}[pink]
				\clip(7,-1.5)rectangle(8.5,1.5);
				\draw (7.8,0) ellipse (2.9cm and 1.15cm);
			\end{scope}
			\begin{scope}[pink]
				\clip(8.5,-1.5)rectangle(14,1.5);
				\draw (9,0) ellipse (2.9cm and 1.13cm);
			\end{scope}
			\begin{scope}[pink]
				\clip(8.5,-1.5)rectangle(14,1.5);
				\draw (9,0) ellipse (2.9cm and 1.15cm);
			\end{scope}

			\draw[brown] (3.3,-1.5) ellipse (0.6cm and 1.5cm);
			\draw[brown] (3.3,-1.5) ellipse (0.6cm and 1.52cm);
			\begin{scope}[magenta]
				\clip(0,-1)rectangle(3.3,-3.5);
				\draw (3.3,-1.5) ellipse (0.65cm and 1.5cm);
			\end{scope}	
			\begin{scope}[magenta]
				\clip(0,-1)rectangle(3.3,-3.5);
				\draw (3.3,-1.5) ellipse (0.65cm and 1.55cm);
			\end{scope}	
			\begin{scope}[magenta]
				\clip(2.65,-1.5)rectangle(3.8,0.7);
				\draw (3,0) ellipse (2cm and 1.1cm);
			\end{scope}	
			\begin{scope}[magenta]
				\clip(2.65,-1.5)rectangle(3.8,0.7);
				\draw (3,0) ellipse (1.9cm and 1.12cm);
			\end{scope}
			\begin{scope}[magenta]
				\clip(2.65,-1.5)rectangle(3.9,0.7);
				\draw (3,0) ellipse (1.9cm and 1.13cm);
			\end{scope}	
			\begin{scope}[magenta]
				\clip(3.3,-5.5)rectangle(8,-1);
				\draw [rotate=-20](4,-.7) ellipse (1.5cm and 1.04cm);
			\end{scope}		
			
			\begin{scope}[magenta]
				\clip(3.3,-5.5)rectangle(8,-1);
				\draw [rotate=-20](4,-.7) ellipse (1.5cm and 1.02cm);
			\end{scope}		
			
			\begin{scope}[lime]
				\clip(3.4,-5.5)rectangle(8,-1.118);
				\draw [rotate=-25](4.9,-.3) ellipse (3cm and 1.04cm);
			\end{scope}	
			\begin{scope}[lime]
				\clip(3.4,-5.5)rectangle(8,-1.118);
				\draw [rotate=-25](4.9,-.3) ellipse (3cm and 1.01cm);
			\end{scope}	
			\begin{scope}[lime]
				\clip(3.3,-1)rectangle(5,-3.5);
				\draw (3.32,-1.5) ellipse (0.65cm and 1.5cm);
			\end{scope}	
			\begin{scope}[lime]
				\clip(3.3,-1)rectangle(5,-3.5);
				\draw (3.32,-1.5) ellipse (0.65cm and 1.51cm);
			\end{scope}
			\begin{scope}[blue]
				\clip(6.4,-5.5)rectangle(15,-1);
				\draw [rotate=-35](7.36,2.45) ellipse (0.8cm and 2.1cm);
			\end{scope}
			\begin{scope}[blue]
				\clip(6.4,-5.5)rectangle(15,-1);
				\draw [rotate=-35](7.35,2.45) ellipse (0.8cm and 2.09cm);
			\end{scope}
			
			\begin{scope}[blue]
				\clip(6.9,-0.97)rectangle(11,1.07);
				\draw (7.425,0) ellipse (2.9cm and 1.092cm);
			\end{scope}
			\begin{scope}[blue]
				\clip(6.9,-0.99)rectangle(11,1.07);
				\draw (7.425,0) ellipse (2.93cm and 1.091cm);
			\end{scope}
			\begin{scope}[blue]
				\clip(7.3,-1.5)rectangle(10,1.1);
				\draw (6.35,0.03) ellipse (2.9cm and 1.16cm);
			\end{scope}
			\begin{scope}[blue]
				\clip(7.3,-1.5)rectangle(10,1.1);
				\draw (6.35,0.03) ellipse (2.9cm and 1.17cm);
			\end{scope}
			\begin{scope}[blue]
				\clip(3.4,-5.5)rectangle(6.5,-1.118);
				\draw [rotate=-25](5,-.3) ellipse (3cm and 1.09cm);
			\end{scope}	
			\begin{scope}[blue]
				\clip(3.4,-5.5)rectangle(6.4,-1.118);
				\draw [rotate=-25](5,-.3) ellipse (3cm and 1.08cm);
			\end{scope}	
			\begin{scope}[blue]
				\clip(3.3,-5.5)rectangle(8,-1);
				\draw [rotate=-20](4,-.7) ellipse (1.6cm and 1.06cm);
			\end{scope}	
			\begin{scope}[pink]
				\clip(7,-1.5)rectangle(10,1);
				\draw (6.27,0.03) ellipse (2.9cm and 1.14cm);
			\end{scope}
			\begin{scope}[pink]
				\clip(7,-1.5)rectangle(10,1);
				\draw (6.27,0.03) ellipse (2.9cm and 1.145cm);
			\end{scope}
			\node at (-0.3,0.3) {$u_{1}^3$};
			\node at (0.5,0){$u_{2}^3$};
			\node at (1.8,0) {$u_{3}^3$};
			\node at (3.9,0.4) {$u_{4}^3$};
			
			\node at (3.3,-0.6) {$u_{5}^3$};
			
			\node at (5.7,0.4) {$u_{6}^3$};
			\node at (5.3,-0.6) {$u_{7}^3$};
			
			\node at (6.7,0.4) {$u_{8}^3$};
			\node at (6.5,-0.6){$u_{9}^3$};
			\node at (8,0) {$u_{10}^3$};
			\node at (9.9,0) {$u_{12}^3$};;
			\node at (11,0) {$u_{13}^3$};
			\node at (10.6,0.5){$u_{11}^3$};
			\node at (10.5,-0.6){$u_{14}^3$};

			\node at (3.4,-1.6) {$u_{15}^3$};
			\node at (3.3,-2.6) {$u_{16}^3$};
			
			\node at (4.4,-1.8) {$u_{17}^3$};
			\node at (4.3,-2.6) {$u_{18}^3$};
			
			\node at (5.4,-2.3) {$u_{19}^3$};
			\node at (6,-3.1) {$u_{20}^3$};
			
			\node at (7.4,-2) {$u_{21}^3$};
			
			\draw [olive] (-1,-4.5) rectangle (12,2);
		\draw[->] [black](0.5,0.7) to [out=10,in=35] (0.2,1.5);
		\node [black]  at (-0.1,1.5) {$N_{1}^1$};
			\draw[->] [yellow](3.1,0.7) to [out=10,in=35] (2.8,1.5);
		\node [yellow]  at (2.5,1.5) {$N_{2}^1$};	
			\draw[->] [red](5.1,0.7) to [out=10,in=35] (4.8,1.5);
		\node [red]  at (4.5,1.5) {$N_{3}^1$};	
			\draw[->] [cyan](6.5,0.7) to [out=10,in=35] (6.2,1.5);
		\node [cyan]  at (5.9,1.5) {$N_{4}^1$};
			\draw[->] [green](8.1,0.6) to [out=10,in=35] (7.8,1.5);
		\node [green]  at (7.5,1.5) {$N_{5}^1$};
			\draw[->] [orange](9.5,0.5) to [out=10,in=35] (9.3,1.4);
		\node [orange]  at (9,1.5) {$N_{6}^1$};
			\draw[->] [pink](11.1,0.5) to [out=10,in=35] (10.8,1.3);
		\node [pink]  at (10.5,1.38) {$N_{7}^1$};
			\draw[->] [brown](3,-2.3) to [out=10,in=35] (2.2,-1.8);
		\node [brown]  at (1.8,-1.8) {$N_{1}^2$};	
			\draw[->] [magenta](3.9,-2.8) to [out=50,in=55] (3.8,-3.7);
		\node [magenta]  at (3.5,-3.8) {$N_{2}^2$};
		\draw[->] [lime](5.5,-3.6) to [out=50,in=55] (5.3,-4.3);
		\node [lime]  at (5,-4.2) {$N_{3}^2$};
			\draw[->] [blue](8,-2.2) to [out=50,in=55] (8.7,-2.6);
		\node [blue]  at (8.7,-2.9) {$N_{4}^2$};
		\end{tikzpicture}\\
		 Figure $8$. Venn Diagram of Figure $7$. 			
	\end{center}
		
		In Figure 8, $N^1_1 N^1_2 N^1_3 N^1_4 N^1_5 N^1_6 N^1_7$ and $N^2_1 N^2_2 N^2_3 N^2_4$ are two nbh chains in $G$. Also,
		\begin{align}
			N^1_2 \setminus N^1_1 &= \{u_{4}^3,u_{5}^3\} \subset \{u_{4}^3,u_{5}^3,u_{6}^3,u_{7}^3\} = N^1_3, \label{26}\\
			N^1_3 \setminus N^1_2 &= \{u_{6}^3,u_{7}^3\} \subset \{u_{6}^3,u_{7}^3,u_{8}^3,u_{9}^3\} = N^1_4,\label{27}\\
			N^1_4 \setminus N^1_3 &= \{u_{8}^3,u_{9}^3\} \subset \{u_{8}^3,u_{9}^3,u_{10}^3\} = N^1_5, \label{28}\\
			N^1_5 \setminus N^1_4 &= \{u_{10}^3\} \subset \{u_{10}^3,u_{12}^3\} = N^1_6, \label{29}\\
			N^1_6 \setminus N^1_5 &= \{u_{12}^3\} \subset \{u_{11}^3,u_{12}^3,u_{13}^3,u_{14}^3\} = N^1_7, \label{34}\\
			N^2_2 \setminus N^2_1 &= \{u_{17}^3,u_{18}^3\} \subset \{u_{17}^3,u_{18}^3,u_{19}^3,u_{20}^3\} = N^2_3, \label{30}\\
			N^2_3 \setminus N^2_2 &= \{u_{19}^3,u_{20}^3\} \subset \{u_{12}^3,u_{19}^3,u_{20}^3,u_{21}^3\} =  N^2_4, \label{31}\\
			N^2_1 \setminus N^2_2 &= \{u_{5}^3\} \subset  \{u_{4}^3,u_{5}^3\} = N^1_3 \setminus N^1_4,\label{32}\\
			N^1_6 \setminus N^1_5 &= \{u_{12}^3\} \subset  \{u_{12}^3,u_{21}^3\} = N^2_4 \setminus N^2_3.\label{33}
		\end{align}
	It is clear that the nbh chains $NC_1$ and $NC_2$ satisfy property (i) of NC-T3 by (\ref{26}) - (\ref{31}). Also, they satisfy property (ii) of NC-T3 by (\ref{32}) and (\ref{33}) for $ i = 3 $, $ j = 1 $  and $ r = 1 $ as given in Definition \ref{d3.1}. Thus, $ NC_1 $ and $ NC_2 $ form an NC-T3.} 
		\end{exm}
	
\begin{theorem}\quad \label{t3} {\rm Let $N^1_1 N^1_2 \cdots N^1_{2n+1}$ and $N^2_1 N^2_2 \cdots N^2_{2n+1}$ be two NC-T1 and form an NC-T2 in a graph $G$. Then the NC-T2 is also an NC-T3 in $G$. 
		
	i.e., An NC-T2 formed on two NC-T1, each of odd length $2n+1$, is also an NC-T3.}
\end{theorem}
\begin{proof}
	Given that $N^1_1 N^1_2 \cdots N^1_{2n+1}$ and $N^2_1 N^2_2 \cdots N^2_{2n+1}$ are NC-T1 and form an NC-T2 in $G$. By  definition, property (ii) of NC-T1 is same as property (i) of NC-T3.
	\begin{center}
		\begin{tikzpicture}[thick,
			set/.style = { circle, minimum size = .061cm}]
			
			\draw[green] (3,0) ellipse (2cm and 1cm);
			\draw[blue] (1,0) ellipse (1.94cm and 0.99cm);
			
			\begin{scope}[red]
				\clip(2,-1)rectangle(3,1.2);
				\draw (1,0) circle [x radius=2cm, y radius=10mm];
			\end{scope}
			\begin{scope}[red]
				\clip(2,-1)rectangle(3,1.2);
				\draw (1,0) circle [x radius=1.975cm, y radius=10mm];
			\end{scope}
			\begin{scope}[blue]
				\clip(2,-1)rectangle(3,1.2);
				\draw (1,0) circle [x radius=1.935cm, y radius=9.8mm];
			\end{scope}

			\begin{scope}[red]
				\clip(2,-1.5)rectangle(7,1.5);
				\draw (4,0) circle [x radius=3cm, y radius=11.9mm];
			\end{scope}
			
			\node[blue] at (.75,-0.6) {$N^1_{i}$};
			\node[green] at (3.5,-0.5) {$N^1_{i+1}$};
			\node[red] at (5.75,0.5) {$N^1_{i+2}$};
			
			\node[blue] at (0,0) {$N^1_{i}\setminus N^1_{i+1}$};
			\node at (2,0) {\tiny $N^1_{i}\cap N^1_{i+1}$};
			\node[green] at (4,0.2) {$N^1_{i+1}\setminus N^1_{i}$};
			\node[red] at (6,-0.2) {$N^1_{i+2}\setminus N^1_{i+1}$};	
		\end{tikzpicture}\\
		\vspace{0.1cm} Figure 9. 			
	\end{center}
	
	In Figure 9, we consider a Venn diagram corresponding to nbhs $N^1_i$, $N^1_{i+1}$, $N^1_{i+2}$ which are represented by regions enclosed by a curve with colors blue, green, red, respectively. We also present property (i) of NC-T3, namely, $ N^1_{i+1} \setminus N^1_{i} \subset N^1_{i+2}$ in the figure for $i$ = 1 to $2n-1$. 
	
	For $i = 1,2,\ldots,2n-1$, we have 
	\begin{align}
		(N^1_{i+1} \setminus N^1_{i}) \subset N^1_{i+2}~ which~ implies, \notag
	\end{align}
 	\begin{align}
	 N^1_{i+1} \setminus N^1_{i+2}& = N^1_{i} \cap N^1_{i+1} ~~and~~\label{aa1}\\ N^1_{i+1} \setminus N^1_{i}& = N^1_{i+1} \cap N^1_{i+2}. ~~~See ~Figure~6.\label{a2}
	\end{align}
		By property (i) of NC-T2,
	\begin{align}
		N^1_{1} \cap N^1_{2} \cap N^2_{1} &= N^2_{1} \setminus N^2_{2} \neq \emptyset ~~~and~\label{a1} \\ N^1_{2n+1} \cap N^2_{2n+1} \cap N^2_{2n} &= N^1_{2n+1} \setminus N^1_{2n} \neq \emptyset. \label{a3}
	\end{align} 
	Also, we have,
	\begin{align}
	N^1_{1} \cap N^1_{2} \cap N^2_{1} \subset N^1_{1} \cap N^1_{2} =  N^1_{2} \setminus N^1_{3}~~~\label{a4} using ~(\ref{aa1}); 
	\end{align} 
	\begin{align}
		N^1_{2n+1} \cap N^2_{2n+1} \cap N^2_{2n} \subset N^2_{2n+1} \cap N^2_{2n} = N^2_{2n} \setminus N^2_{2n-1} ~~~~using ~(\ref{a2}). \label{a5}
	\end{align} 
	From (\ref{a1}),(\ref{a3}),(\ref{a4}) and (\ref{a5}), we get,
	\begin{align}
		\emptyset \neq N^2_{1} \setminus N^2_{2} \subset   N^1_{2} \setminus N^1_{3}~and~~~\emptyset \neq  N^1_{2n+1} \setminus N^1_{2n} \subset N^2_{2n} \setminus N^2_{2n-1}. \label{a6} 
	\end{align} 
$\therefore$ (\ref{a6}) gives the property (ii) of NC-T3, where $i$ = 2, $j$ = 1 and $r = n-1$. Hence  $N^1_1 N^1_2 \cdots N^1_{2n+1}$ and $N^2_1 N^2_2 \cdots N^2_{2n+1}$ form an NC-T3.
\end{proof}

	\begin{theorem}\quad \label{th3} {\rm Any graph $ G $ with an NC-T3 is NDM.}
	\end{theorem}
	
	\begin{proof} \quad
	We prove the theorem by the method of contradiction. Let $N^1_1 N^1_2 \cdots N^1_{n}$ and $N^2_1 N^2_2 \cdots N^2_{p}$ be nbh chains form an NC-T3 in $G$. Here, we consider $n_1 = n$ and  $n_2 = p$.
		
		If $G$ is not an NDM graph, then assume that $G$ be a DM graph. Let $f$ be a DML of $G$ with DM constant $S$. 
		Consider the nbh chain $N^1_1 N^1_2 \cdots N^1_{n}$ in $ G $. By the definition of DML, we get,
		\begin{align}
			S~=\sum_{v_j \in  N^1_1}f(v_j)=~\sum_{v_j \in N^1_2 }f(v_j)~=~\cdots~=~\sum_{v_j \in  N^1_i}f(v_j)~=\notag \\ \cdots~=~\sum_{v_j \in N^1_{n}}f(v_j). \label{1}
		\end{align}
		And for $i = 1,2,\ldots,n-2$,  
		\begin{align}
			N^1_{i+1} \setminus N^1_{i} \subset N^1_{i+2}. \implies N^1_{i+1} \setminus N^1_{i} = N^1_{i+1} \cap N^1_{i+2}.~~See ~Figure~9.\label{2}
		\end{align}
		For $i = 1,2,\ldots,n-1$, we have, 
		\begin{align}
			N^1_{i}&= (N^1_{i} \setminus N^1_{i+1}) \cup (N^1_{i}\cap N^1_{i+1})~and \label{3} \\ 
			N^1_{i+1}&= (N^1_{i+1} \setminus N^1_{i}) \cup (N^1_{i}\cap N^1_{i+1}).\label{4}
		\end{align}
		Also for $i = 1,2,\ldots,n-2$, we have,
		\begin{align}
			N^1_{i+2}&= (N^1_{i+2} \setminus N^1_{i+1}) \cup (N^1_{i+1}\cap N^1_{i+2}) \notag\\
			&= (N^1_{i+2} \setminus N^1_{i+1}) \cup (N^1_{i+1} \setminus N^1_{i}) ~~~ using ~~(\ref{2}).\label{5}
		\end{align}
		\begin{align}
			\Rightarrow	N^1_{i+3} = (N^1_{i+3}\setminus N^1_{i+2}) \cup (N^1_{i+2} \setminus N^1_{i+1}),~ i = 1,2,\dots ,n-3; ~ \label{6}
		\end{align}
		\begin{align}
			N^1_{i+4} = (N^1_{i+4}\setminus N^1_{i+3}) \cup (N^1_{i+3} \setminus N^1_{i+2}),~ i = 1,2,\dots ,n-4 ~ and\label{7} 
		\end{align}
		\begin{align}
			N^1_{i+5} = (N^1_{i+5}\setminus N^1_{i+4}) \cup (N^1_{i+4} \setminus N^1_{i+3}),~ i = 1,2,\dots ,n-5.\label{8} 
		\end{align}
		
		Using (\ref{3})~\&~(\ref{4})~ in  (\ref{1}), for $ i= 1,2,\dots ,n-1 $, we get, \small
		\begin{align}
			S = \sum_{u \in N^1_{i} \setminus N^1_{i+1}}{f(u)}+\sum_{u\in N^1_{i}\cap N^1_{i+1}}{f(u)}
			=\sum_{u\in N^1_{i+1} \setminus N^1_{i}}{f(u)}+\sum_{u\in N^1_{i}\cap N^1_{i+1}}{f(u)}.\notag 
		\end{align}
		\begin{align}
			\Rightarrow	\sum_{v_j\in N^1_{i} \setminus N^1_{i+1}}{f(v_j)}=\sum_{v_j\in N^1_{i+1} \setminus N^1_{i}}{f(v_j)}~~for~~ i= 1,2,\dots ,n-1.\label{9}
		\end{align}
		Using (\ref{5})~\&~(\ref{6})~ in  (\ref{1}), for $ i= 1,2,\dots ,n-3 $, we get,
		\begin{align}
			S   =\sum_{v_j\in N^1_{i+2} \setminus N^1_{i+1}}{f(v_j)} ~+\sum_{v_j\in N^1_{i+1} \setminus N^1_{i}}{f(v_j)} \notag \end{align}
			\begin{align}
			 \hspace{3.5cm} =\sum_{v_j\in N^1_{i+3} \setminus N^1_{i+2}}{f(v_j)}+\sum_{v_j\in N^1_{i+2} \setminus N^1_{i+1}}{f(v_j)}. \notag
		\end{align}
		\begin{align}
			\Rightarrow \sum_{v_j\in N^1_{i+1} \setminus N^1_{i}}{f(v_j)}=\sum_{v_j\in N^1_{i+3} \setminus N^1_{i+2}}{f(v_j)}~~ for ~~i= 1,2,\dots ,n-3. \label{10}
		\end{align}
		Using (\ref{7})~\&~(\ref{8})~ in  (\ref{1}), for $ i= 1,2,\dots ,n-5 $, we get,
		\begin{align}
				S   =\sum_{v_j\in N^1_{i+4} \setminus N^1_{i+3}}{f(v_j)} ~+\sum_{v_j\in N^1_{i+3} \setminus N^1_{i+2}}{f(v_j)},
				 \notag \\=\sum_{v_j\in N^1_{i+5} \setminus N^1_{i+4}}{f(v_j)}+ \sum_{v_j\in N^1_{i+4} \setminus N^1_{i+3}}{f(v_j)}. \notag
		\end{align}
		\begin{align}
			\Rightarrow 	\sum_{v_j\in N^1_{i+3} \setminus N^1_{i+2}}{f(v_j)}=\sum_{v_j\in N^1_{i+5} \setminus N^1_{i+4}}{f(v_j)} ~ for ~ i= 1,2,\dots ,n-5. \label{11}
		\end{align}
		From (\ref{9}), (\ref{10}), (\ref{11}), for $i= 1,2,\dots ,n-5$, we have, 
		\begin{align}
			\sum_{v_j\in N^1_{i} \setminus N^1_{i+1}}{f(v_j)}=\sum_{v_j\in N^1_{i+1} \setminus N^1_{i}}{f(v_j)}=\sum_{v_j\in N^1_{i+3} \setminus N^1_{i+2}}{f(v_j)} =\sum_{v_j\in N^1_{i+5} \setminus N^1_{i+4}}{f(v_j)}.
		\end{align}
		In general, for $ i = 1,2,\ldots, n-1 $ and $ r = 0,1,\ldots,  \left\lfloor{\frac{n}{2}}\right\rfloor -1 
		$, we have,
		\begin{align}
			\sum_{v_j\in N^1_{i} \setminus N^1_{i+1}}{f(v_j)}=\sum_{v_j\in N^1_{i+1+2r} \setminus N^1_{i+2r}}{f(v_j)}. \label{13}
		\end{align}
		Also, for $ j = 1,2,\ldots, p-1 $ and $ r = 0,1,\ldots,  \left\lfloor{\frac{p}{2}}\right\rfloor -1 
		$, corresponding to the neighbourhood chain $N^2_1 N^2_2 \cdots N^2_{p}$ in $ G $, we have,
		\begin{align}
			\sum_{v_j\in N^2_{j} \setminus N^2_{j+1}}{f(v_j)}=\sum_{v_j\in N^2_{j+1+2r} \setminus N^2_{j+2r}}{f(v_j)}. \label{14}
		\end{align}
		 For $ 1 \leq i \leq n-1 $, $ 1 \leq j \leq p-1 $, $ 0 \leq r \leq \left\lfloor{\frac{q}{2}}\right\rfloor -1  $ and  $ q = min ~\{n,p\} $, by property (ii) of NC-T3, and  for a particular set of values of $ i, j$ and $r$, we have, 
		\begin{align}
			N^2_{j} \setminus N^2_{j+1} \subset N^1_{i} \setminus N^1_{i+1} \notag ~and~~
			N^1_{i+1+2r} \setminus 
			N^1_{i+2r} \subset  N^2_{j+1+2r} \setminus N^2_{j+2r}.   \notag
		\end{align}
		\begin{align}
		\implies~	N^1_{i} \setminus N^1_{i+1} &= ((N^1_{i} \setminus N^1_{i+1})\setminus (N^2_{j} \setminus N^2_{j+1})) \cup (N^2_{j} \setminus N^2_{j+1}) ~and \label{17}
		\end{align}
		\begin{align}
			N^2_{j+1+2r} \setminus N^2_{j+2r}&= (( N^2_{j+1+2r} \setminus N^2_{j+2r}) \setminus (N^1_{i+1+2r} \setminus 
			N^1_{i+2r})) \cup ( N^1_{i+1+2r} \setminus 
			N^1_{i+2r}). \label{18}
		\end{align}
		Using (\ref{17}) in (\ref{13}) and (\ref{18}) in (\ref{14}) for $ 1 \leq i \leq n-1 $, $ 1 \leq j \leq p-1 $, $ 0 \leq r \leq \left\lfloor{\frac{q}{2}}\right\rfloor -1  $ and  $ q = min ~\{n,p\} $, we get, 
		\begin{align}
			\sum_{v_j\in (N^1_{i} \setminus N^1_{i+1})\setminus (N^2_{j} \setminus N^2_{j+1})}{f(v_j)} + \sum_{v_j\in N^2_{j} \setminus N^2_{j+1}}{f(v_j)}=\sum_{v_j\in N^1_{i+1+2r} \setminus N^1_{i+2r}}{f(v_j)}; \label{19} 
		\end{align}
		\begin{align}
			\sum_{v_j\in N^2_{j} \setminus N^2_{j+1}}{f(v_j)}=\sum_{v_j\in ( N^2_{j+1+2r} \setminus N^2_{j+2r}) \setminus (N^1_{i+1+2r} \setminus 
				N^1_{i+2r})}{f(v_j)} +\notag \\ \sum_{v_j\in N^1_{i+1+2r} \setminus N^1_{i+2r}}{f(v_j)} \label{20}. 
		\end{align}
			For $ 1 \leq i \leq n-1 $, $ 1 \leq j \leq p-1 $, $ 0 \leq r \leq \left\lfloor{\frac{q}{2}}\right\rfloor -1  $ and  $ q = min ~\{n,p\} $, (\ref{19}) and (\ref{20}) implies,
		\begin{align}
			\sum_{v_j\in (N^1_{i} \setminus N^1_{i+1})\setminus (N^2_{j} \setminus N^2_{j+1})}{f(v_j)} + \sum_{v_j\in ( N^2_{j+1+2r} \setminus N^2_{j+2r}) \setminus (N^1_{i+1+2r} \setminus 
				N^1_{i+2r})}{f(v_j)} = 0, \notag	\end{align}
				which is a contradiction to condition (ii) of NC-T3, since $ ((N^1_{i} \setminus N^1_{i+1})\setminus (N^2_{j} \setminus N^2_{j+1})) \cup (( N^2_{j+1+2r} \setminus N^2_{j+2r}) \setminus (N^1_{i+1+2r} \setminus 
				N^1_{i+2r}))$  contains atleast one element whose labeling is $ \geq 1$.\\
				Hence we get the result. 
	\end{proof}

\begin{theorem}\quad \label{cor1} {\rm Let $G$ be a graph containing two Type-1 nbh chains, each of length $n$, and the two form a nbh chains of Type-2, $n \geq 2$ and $n\in\mathbb{N}$. Then $G$ is NDM.}
\end{theorem}
\begin{proof}
Let $G$ be a graph containing two Type-1 nbh chains, each of length $n$, and the two form a nbh chains of Type-2, $n \geq 2$ and $n\in\mathbb{N}$.

When $n$ is even, $G$ is NDM by  Theorem \ref{t2}. When $n$ is odd, then NC-T2 also forms an NC-T3 in $G$, Theorem \ref{t3}. In this case $G$ is NDM by Theorem \ref{th3}.  
\end{proof}

\begin{theorem} \label{cor2} {\rm For $m \geq 2$ and $n \geq 3$, graph $ P_m \Box C_n $ is NDM, $ m,n \in \mathbb{N} $.}
\end{theorem}
\begin{proof}
	For $m \geq 2$ and $n \geq 3$, using Theorem \ref{t1}, graph $P_{m} \Box C_{n}$ contains a nbh chains of Type-2, $m,n \in \mathbb{N}$. Also, using Theorem \ref{cor1}, $P_{m} \Box C_{n}$ is NDM for $m \geq 2$, $n \geq 3$ and $m,n \in \mathbb{N}$. 
\end{proof}

\begin{note}\quad {\rm Thus, by the above results, Conjecture \ref{c1} is settled completely}.
\end{note}

	\section{Families of NDM graphs by existence of NC-T3}

 Here, we start with a few definitions which are required in this section to obtain families of NDM graphs.

\begin{definition}\cite{xw} For $m\in\mathbb{N}$ and $m \geq 2$,
	let $\{G_i\}^m_{i=1}$ be a collection of graphs with $u_i\in V(G_i)$ as a fixed vertex. \em{The vertex amalgamation}, denoted
	by $Amal(G_i, \{u_i\}, m)$, is a graph formed by taking all the $G_i$'s and identifying $u_i$'s. 
\end{definition}

\begin{definition}\cite{ep92} A \emph{two terminal graph} is a graph with two distinguished vertices, $s$ and $t$ called \emph{source} and \emph{sink}, respectively. 
\end{definition}

\begin{definition}\cite{ep92} \emph{The series composition of 2 two terminal graphs} $X$ and $Y$ is a two terminal graph created from the disjoint union of graphs $X$ and $Y$ by merging the sink of $X$ with the source of $Y$. The source of $X$  becomes the source of series composition and the sink of $Y$ becomes the sink of series composition. 
\end{definition}

\begin{definition}\cite{ep92} \emph{The parallel composition of two terminal graphs} $X$ and $Y$ is a two terminal graph created from the disjoint union of graphs $X$ and $Y$ by merging the sources of $X$ and $Y$ to create the source of parallel composition and merging the sinks of $X$ and $Y$ to create the sink of parallel composition. 
\end{definition}

\begin{definition}\cite{go78} Let $G$ be a graph of order $ n(G)$ and $H$ be a graph with root vertex $v$. Then, the \emph{rooted product graph} of $G$ and $H$ is defined as the graph obtained from $G$ and $H$ by taking one copy of $G$ and $n(G)$ copies of $H$ and identifying the $i^{th}$ vertex of $G$  with the root vertex $v$ in the $i^{th}$ copy of $H$ for every $i\in \{ 1,2,3,...,n(G) \}$.
\end{definition}

 Now, we present the following families of graphs as NDM by the existence of NC-T3 in each of these graphs.

\begin{enumerate}
	\item	Let $ G_1,G_2,\cdots,G_k $ be mutually disjoint graphs such that at least one of them contain an NC-T3. Then $\displaystyle$ $ \bigcup_{i=1}^k G_i $ is NDM.
	
	\item Let $m\in\mathbb{N}$, $m \geq 2$ and $\{G_i\}^m_{i=1}$ be a collection of graphs with $u_i\in V(G_i)$ as a fixed vertex such that at least one of $G_i$ contains an NC-T3. Let $H$ = $Amal(G_i, \{u_i\}, m)$, the vertex amalgamation of $ G_1,G_2,\cdots,G_m $ by taking all the $G_i$'s and identifying $u_i$'s without disturbing the structure of at least one NC-T3 in $H$. Then $H$ is NDM.
		
	\item Let $G$ be a simple graph having an NC-T3, say $NC_1$ and $H$ be a super graph of $G$ without disturbing the structure of $NC_1$ in $H$. Then $H$ is NDM.
	
	\item Let $G$ be a simple graph containing an NC-T3, say $NC_1$ and $H$ be a super graph of $G$ containing $NC_1$ as its induced subgraph. Then $H$ is NDM.
	
	\item	Let $ G $ be a simple graph containing an NC-T3, say $NC_1$. Then any subgraph $H$ of $G$ containing $NC_1$ is also NDM.
	
	\item	Let $G$ be a two terminal simple graph and $H$ be a graph having NC-T3. Then the series composition of $G$ and $H$ is NDM, follows from the structure of the combined graph.
	
	\item	Let $G$ be a two terminal simple graph and $H$ be a graph having NC-T3. Then the parallel composition of $G$ and $H$ is NDM, follows from the structure of the combined graph.
	
	\item	Let $G$ be any simple graph and $H$ be a graph having NC-T3. Then the rooted product of $G$ and $H$ is NDM.
	
\end{enumerate}

\noindent
{\bf Conclusion.} \quad In this paper, we obtained families of non-distance magic graphs by the existence of NC-T3 in these graphs and also completely settled the conjecture on NDM of cylindrical grid graph $P_m \Box C_n$ for $m \geq 2$, $n \geq 3$ and $m,n \in \mathbb{N}$.  

\vspace{.5cm}
\noindent
\textbf{Declaration of competing interest}  The authors declare no competing interests.

\vspace{.5cm}
\noindent
\textbf{Acknowledgement} \textit {\rm We express our sincere thanks to the Central University of Kerala, Kasaragod, Kerala, India for providing facilities to carry out this research work and the first author acknowledges the support by Council of Scientific and Industrial Research (CSIR), Government of India through its SRF fellowship.}

\end{document}